\documentclass{article}
\usepackage[utf8]{inputenc}
\usepackage[T1]{fontenc}
\usepackage{graphicx}
\usepackage{amsmath,amssymb,amsfonts}%
\usepackage{a4wide}
\usepackage[round]{natbib}
\usepackage{authblk}
\usepackage{hyperref}

\title{The distributions of convex body functionals of the Poisson~polyhedron}
\author[1]{Felix Ballani\footnote{Email address: \texttt{f.ballani@hzdr.de} | ORCID-iD: \href{https://orcid.org/0000-0002-0415-0220}{0000-0002-0415-0220}}}
\author[2]{Dietrich Stoyan\footnote{Email address: \texttt{stoyan@math.tu-freiberg.de}}}
\affil[1]{Helmholtz-Zentrum Dresden -- Rossendorf, Helmholtz Institute Freiberg for Resource Technology, Chemnitzer Str. 40, 09599 Freiberg, Germany}
\affil[2]{Institut f\"ur Stochastik, TU Bergakademie Freiberg, 09596 Freiberg,
Germany}
\date{}
\begin{document}

\maketitle

\begin{abstract}
This note presents new findings on the distribution of the volume and surface area of both the Poisson polyhedron and the Crofton polyhedron. It also considers the distribution of the area of the Poisson polygon and the Crofton polygon in the plane. This note motivates the modeling of the distributions of these convex body functionals as mixture distributions and the essential inclusion of the generalized gamma distribution therein as mixture components.
Finally, model parameters are estimated and models are selected based on extensive stochastic simulations of these random polytopes. This includes that for the first time, this note gives very accurate estimates of the side number probabilities for both the Poisson polyhedron and the Crofton polyhedron, as these are used as mixing weights in the distribution models.
\end{abstract}

\section{Introduction}\label{sec:intro}
The Poisson polyhedron is an important model for random polyhedra. It is used in the context of crushing hard materials, which in several cases is modeled as tessellating the space by random planes. The resulting fragments are then considered the generated particles \citep{MottLinfoot1943,Grady2009,Grady2017,BallaniStoyan2025}. For these particles, volumes and surface areas are of interest. Conversely, the Poisson polyhedron itself is a building block in various modeling approaches. For instance, it serves as a particle model in the context of investigations of two-phase materials \citep{Serra1982,QuenecEtAl1994,CosterChermant2002,CSKM2013,EscodaEtAl2015,EscodaEtAl2016,Jeulin2021}, and it is also used as a basis for extended particle models \citep{BallaniBoogaart2014}. Furthermore, the Poisson polyhedron is an important part of the construction of certain max-stable random fields \citep{LantuejoulEtAl2011}.

The Poisson polyhedron appears in the context of the well-known Poisson plane tessellation, which is generated by a Poisson plane process. This process divides space into an infinite number of convex cells. If the plane process is stationary, it makes sense to speak of ``the typical cell'', a random polyhedron corresponding to a random selection principle in which a cell is randomly selected from the cells of the entire tessellation, with each cell having the same chance of being selected. (For precise definitions, see \citet{SchneiderWeil2008}, \citet{CSKM2013}, or \citet{HugSchneider2024}.) If the plane process is also isotropic, the polyhedron is called ``Poisson polyhedron'' or ``Poisson cell''. It is important to note that the Poisson polyhedron is also the typical cell of a broader class of tessellations known as STIT-tessellations \citep{NagelWeiss2003,NagelWeiss2005}. These tessellations may be even more realistic in the context of comminution because they model time-dependent, stepwise subdivisions of space.
 
The Poisson polyhedron has random volume, surface area, and mean width. The aim of this note is to determine the distributions of these convex body functionals for the spatial case as introduced here and for the planar case as well. It closes a gap in the stochastic geometrical literature: While there are formulas for some moments and even results on the distribution of the mean width (see, for example, the recent monograph \citet{HugSchneider2024}), it lacks information on the distributions of volume and surface area, although there are algorithms for simulating the typical cell. 

\section{Material and Methods}
\subsection{The model}
Let $C$ be the typical cell (denoted as Poisson cell) and $C_o$ the zero cell (denoted as Crofton cell) of a stationary and isotropic Poisson hyperplane tessellation in $\mathbb{R}^d$, ${d=2,3}$, which is generated by a stationary and isotropic Poisson hyperplane process of intensity ${\lambda>0}$. The parameter $\lambda$ is chosen so that it corresponds to the mean total length of the underlying line process per unit area (${d=2}$) or the mean total area of the underlying plane process per unit volume (${d=3}$) \citep[see, e.\,g.,][Chapters 8 and 9]{CSKM2013}.

In the following, for a convex polygon $K$, $A(K)$ denotes its area, $L(K)$ its perimeter,\linebreak ${B(K)=L(K)/\pi}$ its mean width and $N(K)$ the number of its edges. Similarly, for a convex polyhedron $K$, $V(K)$ is its volume, $S(K)$ its surface area, $B(K)$ its mean width and $N(K)$ the number of its side faces.

\subsubsection{Description and known results}
The following uses the generalized gamma distribution for modeling the distributions of volume and surface area of the Poisson polyhedron. The generalized gamma distribution $\mathsf{GG}(\beta,\nu,\kappa)$ goes back to \citet{Stacy1962} and is given by the density function 
\begin{equation}\label{eq:GG:density}
f(x\mid\beta,\nu,\kappa)=\frac{\kappa}{\beta^{\nu}\Gamma\left(\frac{\nu}{\kappa}\right)}\,x^{\nu-1}\mathrm{e}^{-\left(\frac{x}{\beta}\right)^{\kappa}},\quad x>0
\end{equation} 
with scale parameter ${\beta>0}$ and two shape parameters ${\nu>0}$ and ${\kappa>0}$. It includes many common models of nonnegative random variables, such as the gamma distribution (${\kappa=1}$) and the Weibull distribution (${\nu=\kappa}$). It is therefore a flexible model and thus generally a good candidate for use in modeling the volume and surface area distributions of the Poisson polyhedron.

The following considerations show the specific way in which the generalized gamma distribution is incorporated into the modeling.
\begin{enumerate}
\item The starting point is the finding of \citet{Miles1971} and \citet{MollerZuyev1996} that (in all dimensions $d$) the mean width $B(C)$ of the Poisson cell $C$ with a given number of sides ${N(C)=n}$ follows a gamma distribution ${\mathsf{Gamma}(n-d,\lambda)}$ with shape parameter ${n-d}$ and rate parameter $\lambda$. 
\item Furthermore, if $X$ is a random variable with a gamma distribution with shape parameter $\alpha$ and rate parameter $\lambda$, then $X^k$, for ${k>0}$, has a ${\mathsf{GG}(\lambda^{-k}, \alpha/k, 1/k)}$-distribution. 
\item The volume of $C$, $V(C)$, can be written as the product ${V(C)=Z(C)B(C)^d}$ with the shape factor ${Z(C)=V(C)/B(C)^d}$, which depends only on the shape of $C$, not its size \citep{BonnetEtAl2018,HugSchneider2024}.
\item Furthermore, under the condition ${N(C)=n}$ the mean width $B(C)$ and the shape of $C$ are stochastically independent \citep{MollerZuyev1996,BonnetEtAl2018}. 
\item In summary, the conditional random variable ${V(C)\,\vert\,N(C)=n}$ can be represented as a product of the ${\mathsf{GG}(\lambda^{-d}, (n-d)/d, 1/d)}$-distributed quantity ${B(C)^d\,\vert\,N(C)=n}$ and an independent shape factor ${Z(C)\,\vert\,N(C)=n}$. 
\end{enumerate}
This shows that ${V(C)\,\vert\,N(C)=n}$ generally does not follow a generalized gamma distribution. Nevertheless, the generalized gamma distribution is a natural candidate for approximation. The shape factor $Z(C)$ is positive and, due to the isoperimetric inequality, cannot exceed $\pi/4$ (${d=2}$) or $\pi/6$ (${d=3}$). As the number of sides $n$ increases, the simulations indicate that ${Z(C)\,\vert\,N(C)=n}$ is not only larger on average \citep[the cells are rounder, see also the conjecture in][p. 327]{HugSchneider2024}, but also scatters less. Therefore, it can be expected that for the larger $n$ the random influence of the shape factor will decrease, causing the distribution of ${V(C)\,\vert\,N(C)=n}$ to correspond more and more to the distribution of ${B(C)^d\,\vert\,N(C)=n}$ apart from a scaling in the order of magnitude of ${Z(C)\,\vert\,N(C)=n}$.

Inspired by the decomposition of the distribution of the mean width $B(C)$ into a mixture of conditional distributions with respect to the number of sides, $N(C)$, the approach used in the following is to represent the distribution of the volume $V(C)$ also as a mixture of the distributions of the conditional random variables ${V(C)\,\vert\,N(C)=n}$ with mixture weights ${p(n)=P(N(C)=n)}$ and to approximate the components of the mixture by generalized gamma distributions. For the distribution of the area of the Poisson cell in the case ${d=2}$, \citet{Calka2003b} already pursues such a decomposition. However, his formulas do not lead to a practically usable distribution. Therefore, the approach of this note is to simulate a large number of realizations of the Poisson cell and to fit generalized gamma distributions to the observations of ${V(C)\,\vert\,N(C)=n}$.

Knowing the probabilities ${p(n)=P(N(C)=n)}$, ${n=d+1,d+2,\ldots}$, is therefore necessary for modeling the overall distribution. Since these probabilities do not seem to follow any known parametric model, they are determined by simulation. Their values for ${d=2}$ were determined by \citet{CrainMiles1976,George1987,MichelParoux2007}. Furthermore, exact values are known for side numbers 3 (${p(3)=2-\pi^2/6\approx0.355066}$) \citep{Miles1964a} and 4 (${p(4)\approx0.381466}$) \citep{Tanner1983a}. \citet{Calka2003b} provides a general integral formula for these probabilities and determines them numerically for the first number of sides. However, the values of \citet{CrainMiles1976,George1987,Calka2003b,MichelParoux2007} differ slightly. Therefore, for the case ${d=2}$, separate Monte Carlo estimates of these probabilities $p(n)$ are given here for comparison purposes. For the case ${d=3}$, which is of primary interest here, the authors are unaware of any explicit values for $p(n)$ in the literature; however, see the bar graph in \citet[][Fig. 7a]{KlattEtAl2017}. These values are probably being reported explicitly here for the first time.

The argument for $V(C)$ is transferred to the distribution of the surface area $S(C)$. In particular, for ${d=3}$ the representation ${S(C)=\widetilde{Z}(C)B(C)^2}$ with a shape factor ${\widetilde{Z}(C)=S(C)/B(C)^2}$ applies, where $0\leq\widetilde{Z}(C)\leq\pi$. Therefore, ${S(C)\,\vert\,N(C)=n}$ can be represented as a product of the ${\mathsf{GG}(\lambda^{-2}, (n-3)/2, 1/2)}$-distributed quantity ${B(C)^2\,\vert\,N(C)=n}$ and an independent shape factor ${\widetilde{Z}(C)\,\vert\,N(C)=n}$.

Much the same applies to the Crofton polyhedron $C_o$. The mean width $B(C_o)$ for a given number of sides ${N(C_o)=n}$ follows a gamma distribution with shape parameter $n$ and rate parameter $\lambda$ \citep{MollerZuyev1996}. Therefore, the previous approach regarding the distribution of the Poisson polyhedron's volume can analogously be applied to the Crofton polyhedron. Again, the probabilities ${p_o(n)=P(N(C_o)=n)}$, ${n=d+1,d+2,\ldots}$, do not seem to follow any known parametric model. In the case of ${d=2}$, in addition to the exact value ${p_o(3)\approx0.07682}$ \citep{MichelParoux2007} for the $p_o(n)$ numerically determined values \citep{Calka2003b} and Monte Carlo estimates are known \citep{CrainMiles1976,George1987,MichelParoux2007}. For ${d=3}$, no such information is available.

\subsubsection{Simulation}
For the simulation of the Poisson cell $C$ there exist three different methods. The first generates a realization of a stationary and isotropic Poisson hyperplane tessellation in a large area. Then a tessellation cell is randomly selected, with each cell having an equal chance of being chosen \citep{CrainMiles1976, MichelParoux2007}. 

The second method generates the Poisson cell as a subpolytope of the zero cell $C_o$. It goes back to \citet[p. 221]{Miles1974}, see also \citet[Sect. 12.3.2]{Lantuejoul2002} for a description and \citet[Thm. 10.4.7]{SchneiderWeil2008} for the theoretical background. The standard algorithm for generating the zero cell generates sequentially independent hyperplanes of the Poisson hyperplane process with increasing distance from the origin, where -- in addition to independent uniform directions -- the distances come from a stationary Poisson point process on $[0,\infty)$ of intensity $2\lambda$. This process continues until the farther hyperplanes do not intersect the origin-covering polytope formed by all previously generated hyperplanes. 

The third method \citep{Calka2001,Calka2010} uses the property of the Poisson cell that the diameter of the largest inscribed sphere has an exponential distribution with a mean value $1/\lambda$ \citep[see, e.\,g.,][]{MollerZuyev1996}. Starting from an arbitrary random simplex with such an inscribed sphere centered at the origin, the Poisson cell is the intersection of this simplex and the polytope that contains the origin and results as a half-space intersection induced by those hyperplanes of a stationary and isotropic Poisson hyperplane process of intensity $\lambda$ whose distances to the origin are greater than the radius of the sphere; see also \citet[Thm. 10.4.6]{SchneiderWeil2008} for a justification.

All results in this paper are based on simulations of the Poisson cell according to Calka's method (third method) or on simulations of the zero cell according to the standard method described above.

\subsection{Simulations and their validation}
The simulations were carried out using a program in the Julia language \citep{Ballani2024c}. For both dimensions, 50,000,000 realizations of the Poisson or Crofton cell were generated. The following were determined: the area, boundary length, and number of sides for $d=2$, as well as the volume, surface area, mean width, and number of sides for $d=3$. In all cases, the parameter value $\lambda=1$ was used. The data sets of polygon and polyhedron statistics obtained in this way, on which the following results are based, are published \citep{Ballani2024a} and freely accessible.

Comparisons of theoretically known moments \citep[see, e.\,g.,][]{CSKM2013,MichelParoux2007} with the corresponding empirical moments estimated from the simulations including the estimated standard errors in Tables \ref{tab:moments:d2} and \ref{tab:moments:d3} and goodness-of-fit tests in Tables \ref{tab:KS:d2} and \ref{tab:KS:d3} for the existence of a ${\mathsf{Gamma}(n-d,1)}$-distribution for ${B(C)\,\vert\,N(C)=n}$ or a ${\mathsf{Gamma}(n,1)}$-distribution for ${B(C_o)\,\vert\,N(C_o)=n}$, ${n=d+1,d+2,\ldots}$, suggest that the simulations are reliable and therefore suitable for obtaining further results.

\begin{table}
\centering
\begin{tabular}{r|r|rr}
$E(\ldots)$&theoretical&simulated&s.\,e.\\
\hline
$A(C)$&3.1416&3.1413&0.0009\\
$L(C)$&6.2832&6.2826&0.0008\\
$N(C)$&4.0000&3.9998&0.0001\\
$A(C)^2$&48.7045&48.6850&0.0460\\
$L(C)^2$&68.4438&68.4384&0.0173\\
$N(C)^2$&16.9348&16.9325&0.0012\\
\hline
$A(C_o)$&15.5031&15.5020&0.0025\\
$L(C_o)$&15.5031&15.5031&0.0011\\
$N(C_o)$&4.9348&4.9348&0.0002\\
$A(C_o)^2$&549.3653&549.0751&0.2459
\end{tabular}
\caption{Known theoretical first and second moments (rounded values) of area ($A$), boundary length ($L$) and number of sides ($N$) of the Poisson cell $C$ and the Crofton cell $C_o$ in dimension ${d=2}$ \citep{MichelParoux2007,CSKM2013} and Monte Carlo estimates together with standard errors (s.\,e.).}
\label{tab:moments:d2}
\end{table}

\begin{table}
\centering
\begin{tabular}{r|r|rr}
$E(\ldots)$&theoretical&simulated&s.\,e.\\
\hline
$V(C)$&15.2789&15.285&0.0075\\
$S(C)$&30.5577&30.5639&0.0076\\
$B(C)$&3.0000&3.0002&0.0003\\
$N(C)$&6.0000&6.0001&0.0002\\
$V(C)^2$&3072.0000&3079.1335&8.4466\\
$S(C)^2$&3840.0000&3842.8622&2.9392\\
$B(C)^2$&14.6921&14.7035&0.0095\\
$N(C)^2$&38.6921&38.6940&0.0032\\
\hline
$V(C_o)$&201.0619&200.9230&0.1005\\
$S(C_o)$&201.0619&200.9657&0.0566\\
$B(C_o)$&8.5797&8.5779&0.0011\\
$N(C_o)$&8.5797&8.5783&0.0006
\end{tabular}
\caption{Known theoretical first and second moments (rounded values) of volume ($V$), surface area ($S$), mean width ($B$) and number of sides ($N$) of the Poisson cell $C$ and the Crofton cell $C_o$ in dimension ${d=3}$ \citep{CSKM2013} and Monte Carlo estimates together with standard errors (s.\,e.).}
\label{tab:moments:d3}
\end{table}

\begin{table*}
\centering
\begin{tabular}{r|r|r|r|r|r|r|r|r|r|r}
$n$&3&4&5&6&7&8&9&10&11&12\\
\hline
$C$&0.75&0.93&0.20&0.05&0.34&0.35&0.57&0.03&0.35&0.26\\
\hline
$C_o$&0.17&0.18&0.92&0.08&0.88&0.83&0.25&0.20&0.05&0.57
\end{tabular}
\caption{$p$-values of the Kolmogorov-Smirnov test to test the hypotheses that in dimension ${d=2}$ ${B(C)\,\vert\,N(C)=n}$ follows a ${\mathsf{Gamma}(n-2,1)}$-distribution or ${B(C_o)\,\vert\,N(C_o)=n}$ follows a ${\mathsf{Gamma}(n,1)}$-distribution.}
\label{tab:KS:d2}
\end{table*}

\begin{table*}
\centering
\begin{tabular}{r|r|r|r|r|r|r|r|r|r|r}
$n$&4&5&6&7&8&9&10&11&12&13\\
\hline
$C$&0.48&0.19&0.17&0.89&0.86&0.51&0.72&0.82&0.92&0.68\\
\hline
$C_o$&0.26&0.31&0.16&0.20&0.25&0.72&0.19&0.97&0.69&0.44
\end{tabular}
\begin{tabular}{r|r|r|r|r|r|r|r}
$n$&14&15&16&17&18&19&20\\
\hline
$C$&0.14&0.85&0.45&0.94&0.02&0.51&0.28\\
\hline
$C_o$&0.67&0.63&0.28&0.97&0.16&0.63&0.56
\end{tabular}
\caption{$p$-values of the Kolmogorov-Smirnov test to test the hypotheses that in dimension ${d=3}$ ${B(C)\,\vert\,N(C)=n}$ follows a ${\mathsf{Gamma}(n-3,1)}$-distribution or ${B(C_o)\,\vert\,N(C_o)=n}$ follows a ${\mathsf{Gamma}(n,1)}$-distribution.}
\label{tab:KS:d3}
\end{table*}

\section{Results}
The following presents the results of the Poisson and Crofton cell simulations in dimensions $d=2$ and $d=3$, as well as the subsequent maximum likelihood (ML) estimates of the generalized gamma distribution parameters. These ML estimates are provided in Tables \ref{tab:gg:A:d2}, \ref{tab:gg:A0:d2}, \ref{tab:gg:V:d3}, \ref{tab:gg:S:d3}, \ref{tab:gg:V0:d3} and \ref{tab:gg:S0:d3}. To assess their precision, the asymptotic standard errors resulting from the Fisher information of the generalized gamma distribution \citep{RebbahEtAl2019} are also given. The tables also show the numbers $m(n)$ (Poisson cell) or $m_o(n)$ (Crofton cell) of the realizations with exactly $n$ sides and thus the sample size on which the ML estimates are based.
\subsection{The case $d=2$}
Table \ref{tab:prob:d2} provides Monte Carlo estimates $m(n)/M$ and $m_o(n)/M$ of the probabilities\linebreak ${p(n)=P(N(C)=n)}$ and ${p_o(n)=P(N(C_o)=n)}$ together with the associated empirical standard errors. These values show excellent agreement within the expected precision with the theoretically known \citep{Miles1964a,Tanner1983a,MichelParoux2007} and empirically determined \citep{CrainMiles1976,George1987,MichelParoux2007} values.

Tables \ref{tab:gg:A:d2} and \ref{tab:gg:A0:d2} then provide the ML estimates ${(\hat{\beta}(n),\hat{\nu}(n),\hat{\kappa}(n))}$ and ${(\hat{\beta}_o(n),\hat{\nu}_o(n),\hat{\kappa}_o(n))}$ of the parameters ${(\beta,\nu,\kappa)}$ of the generalized gamma distribution for ${A(C)\,\vert\,N(C)=n}$ and\linebreak ${A(C_o)\,\vert\,N(C_o)=n}$, respectively, together with related standard errors. Furthermore, the $p$-values of a Kolmogorov-Smirnov (KS) test are given for the hypothesis that the area of polygons with a number of sides of $n$ has a ${\mathsf{GG}(\hat{\beta}(n), \hat{\nu}(n),\hat{\kappa}(n))}$-distribution or a ${\mathsf{GG}(\hat{\beta}_o(n),\hat{\nu}_o(n),\hat{\kappa}_o(n))}$-\linebreak distribution. The tests show that the conditional distributions ${A(C)\,\vert\,N(C)=n}$ or\linebreak ${A(C_o)\,\vert\,N(C_o)=n}$ can indeed be modeled sufficiently well by the corresponding generalized gamma distributions. 

\begin{table*}
\centering
\begin{tabular}{r|rr|rr}
$n$&$\hat{p}(n)$&s.\,e.&$\hat{p}_o(n)$&s.\,e.\\
\hline
3&0.355130&0.000068&0.076832&0.000038\\
4&0.381482&0.000069&0.300454&0.000065\\
5&0.189570&0.000055&0.342646&0.000067\\
6&0.058767&0.000033&0.193892&0.000056\\
7&0.012721&0.000016&0.067136&0.000035\\
8&0.002047&0.000006&0.015887&0.000018\\
9&0.000255&0.000002&0.002750&0.000007\\
10&0.000025&0.000001&0.000361&0.000003
\end{tabular}
\caption{Monte Carlo estimates of the probabilities ${p(n)=P(N(C)=n)}$ and ${p_o(n)=P(N(C_o)=n)}$ as well as corresponding standard errors (s.\,e.) in dimension ${d=2}$. Theoretically known values are ${p(3)=2-\pi^2/6=0.3550659\ldots}$ \citep{Miles1964a} and ${p(4)=0.3814662\ldots}$ \citep{Tanner1983a} as well as ${p_o(3)=0.07682\ldots}$ \citep{MichelParoux2007}. Estimated probabilities for the higher side numbers not listed in the table: ${\hat{p}(11)=1.86\text{e-}6}$, ${\hat{p}(12)=1.2\text{e-}7}$, ${\hat{p}_o(11)=3.706\text{e-}5}$, ${\hat{p}_o(12)=3.22\text{e-}6}$, ${\hat{p}_o(13)=2.8\text{e-}7}$; higher side numbers did not occur.}\label{tab:prob:d2}
\end{table*}

\begin{table*}
\centering
\begin{tabular}{r|r|rr|rr|rr|r}
$n$&$m(n)$&$\hat{\beta}(n)$&s.\,e.&$\hat{\nu}(n)$&s.\,e.&$\hat{\kappa}(n)$&s.\,e.&KS test\\
\hline
3&17756516&0.2710&0.0009&0.5029&0.0003&0.4701&0.0004&0.1220\\
4&19074115&0.2950&0.0016&1.0015&0.0007&0.4637&0.0004&0.4598\\
5&9478503&0.3402&0.0035&1.5038&0.0020&0.4685&0.0007&0.9170\\
6&2938325&0.3993&0.0088&1.9964&0.0055&0.4761&0.0015&0.5393\\
7&636073&0.4150&0.0232&2.5036&0.0167&0.4761&0.0035&0.9895\\
8&102367&0.4223&0.0678&3.0317&0.0560&0.4761&0.0095&0.9981\\
9&12735&0.4855&0.2389&3.5202&0.1978&0.4846&0.0291&0.9877\\
10&1267&0.2776&0.5543&4.3856&0.9024&0.4603&0.0997&0.9736
\end{tabular}
\caption{ML estimates of the parameters $\beta$, $\nu$ and $\kappa$ of the generalized gamma distribution for ${A(C)\,\vert\,N(C)=n}$, the corresponding standard errors (s.\,e.) as well as the $p$-values of each a KS test.}\label{tab:gg:A:d2}
\end{table*}

\begin{table*}
\centering
\begin{tabular}{r|r|rr|rr|rr|r}
$n$&$m_o(n)$&$\hat{\beta}_o(n)$&s.\,e.&$\hat{\nu}_o(n)$&s.\,e.&$\hat{\kappa}_o(n)$&s.\,e.&KS test\\
\hline
3&3841598&0.2742&0.0044&1.4974&0.0031&0.4703&0.0011&0.8854\\
4&15022680&0.3162&0.0032&1.9891&0.0024&0.4686&0.0006&0.8350\\
5&17132342&0.3551&0.0039&2.4952&0.0032&0.4714&0.0007&0.3570\\
6&9694596&0.4011&0.0066&2.9914&0.0056&0.4761&0.0010&0.2913\\
7&3356806&0.4297&0.0132&3.4983&0.0121&0.4785&0.0018&0.8828\\
8&794362&0.5275&0.0351&3.9371&0.0296&0.4891&0.0039&0.9252\\
9&137500&0.3767&0.0723&4.5893&0.0920&0.4688&0.0099&0.8838\\
10&18088&0.7782&0.3649&4.7196&0.2534&0.5088&0.0288&0.8881
\end{tabular}
\caption{ML estimates of the parameters $\beta$, $\nu$ and $\kappa$ of the generalized gamma distribution for ${A(C_o)\,\vert\,N(C_o)=n}$, the corresponding standard errors (s.e.) as well as the $p$-values of each a KS test.}\label{tab:gg:A0:d2}
\end{table*}

Due to the flexibility of the generalized gamma distribution, one could imagine modeling the distribution of $A(C)$ using only one generalized gamma distribution. In this case, the ML estimates are ${\hat{\beta}=2.4180}$, ${\hat{\nu}=0.5079}$, and ${\hat{\kappa}=0.5825}$. A comparison with the values in Table \ref{tab:gg:A:d2} suggests that, in particular, the value of ${\hat{\nu}=0.5079}$ was strongly influenced by the areas of the triangular Poisson cells (${n=3}$); otherwise, the necessary polynomial behavior of the density function for arguments near 0 cannot be guaranteed. However, from a formal point of view the KS test with a $p$-value of virtually 0 clearly rejects the fitted generalized gamma distribution. The same observation can be made for the distribution of $A(C_o)$. The ML estimates of the parameters of a single generalized gamma distribution are ${\hat{\beta}_o=1.8342}$, ${\hat{\nu}_o=1.5536}$, ${\hat{\kappa}_o=0.5478}$; the $p$-value of the corresponding KS test is $5.45\text{e-}8$.

Since a simple distribution model is obviously not sufficient, this note finally proposes to use a mixed model with a density of the form
\begin{equation}\label{eq:mod:A}
\begin{aligned}
f_{A(C)}^{(k)}(x)=\sum_{n=3}^{k}\frac{\hat{p}(n)}{c(k)}\,f\bigg(x\;\bigg|\;\frac{\hat{\beta}(n)}{\lambda^2},\hat{\nu}(n),\hat{\kappa}(n)\bigg),&\\
\quad x>0,&
\end{aligned}
\end{equation}
for the distribution of $A(C)$ with the parameter estimates from Tables \ref{tab:prob:d2} and \ref{tab:gg:A:d2} and the sum $c(k)$ of all $\hat{p}(n)$ for $3\leq n\leq k$. Among the models $f_{A(C)}^{(k)}$, ${k=3,4,\ldots}$, $f_{A(C)}^{(9)}$ has the smallest value AIC for the Akaike information criterion. In addition, the KS test yields a $p$-value of about 0.34, so that $f_{A(C)}^{(9)}$ is also an acceptable model from a formal point of view. 

Among all models with a density of the form
\begin{equation}\label{eq:mod:A0}
\begin{aligned}
f_{A(C_o)}^{(k)}(x)=\sum_{n=3}^{k}\frac{\hat{p}_o(n)}{c_o(k)}\,f\bigg(\!x\,\bigg|\,\frac{\hat{\beta}_o(n)}{\lambda^2},\hat{\nu}_o(n),\hat{\kappa}_o(n)\!\!\bigg),&\\
\quad x>0,&    
\end{aligned}
\end{equation}
with the parameter estimates from Tables \ref{tab:prob:d2} and \ref{tab:gg:A0:d2} and the sum $c_o(k)$ of all $\hat{p}_o(n)$ for $3\leq n\leq k$, the model $f_{A(C_o)}^{(10)}$ has the smallest AIC value and is acceptable with a $p$-value of approximately 0.56 in the KS test. 

\subsection{The case $d=3$}
Table \ref{tab:prob:d3} provides Monte Carlo estimates $m(n)/M$ and $m_o(n)/M$ of the probabilities\linebreak ${p(n)=P(N(C)=n)}$ and ${p_o(n)=P(N(C_o)=n)}$ together with the associated empirical standard errors.
\subsubsection{The Poisson polyhedron}
Using only a single generalized gamma distribution to model the distribution of $V(C)$ results in the ML estimates ${\hat{\beta}_V=6.3503}$, ${\hat{\nu}_V=0.3405}$, and ${\hat{\kappa}_V=0.4060}$. However, the KS test clearly rejects this model ($p$-value of virtually 0). The same applies to $S(C)$ (ML estimates ${\hat{\beta}_S=32.1892}$, ${\hat{\nu}_S=0.5226}$, ${\hat{\kappa}_S=0.6760}$).
Therefore, again a more complex distribution model is required in each case.

Tables \ref{tab:gg:V:d3} and \ref{tab:gg:S:d3} then provide the ML estimates ${(\hat{\beta}_V(n),\hat{\nu}_V(n),\hat{\kappa}_V(n))}$ and ${(\hat{\beta}_S(n),\hat{\nu}_S(n),\hat{\kappa}_S(n))}$ of the parameters ${(\beta,\nu,\kappa)}$ of the generalized gamma distribution for ${V(C)\,\vert\,N(C)=n}$ and\linebreak ${S(C)\,\vert\,N(C)=n}$, respectively, for $n=4,\ldots,16$, together with related standard errors. In addition, the $p$-values of a KS test are provided for the hypothesis that the volume of the polyhedra with a number of sides of $n$ has a ${\mathsf{GG}(\hat{\beta}_V(n), \hat{\nu}_V(n),\hat{\kappa}_V(n))}$-distribution and that the surface area of the polyhedra with a number of sides of $n$ has a ${\mathsf{GG}(\hat{\beta}_S(n),\hat{\nu}_S(n),\hat{\kappa}_S(n))}$-distribution. The test results show that the conditional distributions ${V(C)\,\vert\,N(C)=n}$ or ${S(C)\,\vert\,N(C)=n}$ can indeed be modeled sufficiently well by the corresponding generalized gamma distributions.

\medskip

Distribution of $V(C)$: Among the mixed models
\begin{equation}\label{eq:mod:V}
\begin{aligned}
f_{V(C)}^{(k)}(x)=\sum_{n=4}^{k}\frac{\hat{p}(n)}{c(k)}\,f\bigg(x\,\bigg|\,\frac{\hat{\beta}_V(n)}{\lambda^3},\hat{\nu}_V(n),\hat{\kappa}_V(n)\!\bigg),&\\
\quad x>0,&
\end{aligned}
\end{equation}
with the parameter estimates from Tables \ref{tab:prob:d3} and \ref{tab:gg:V:d3} and the sum $c(k)$ of all $\hat{p}(n)$ for $4\leq n\leq k$, the model $f_{V(C)}^{(16)}$ has the smallest AIC value and is also formally accepted with a $p$-value of 0.47. 

\medskip

Distribution of $S(C)$: Among the mixed models
\begin{equation}\label{eq:mod:S}
\begin{aligned}
f_{S(C)}^{(k)}(x)=\sum_{n=4}^{k}\frac{\hat{p}(n)}{c(k)}\,f\bigg(x\;\bigg|\;\frac{\hat{\beta}_S(n)}{\lambda^2},\hat{\nu}_S(n),\hat{\kappa}_S(n)\bigg),&\\
\quad x>0,&
\end{aligned}
\end{equation}
with $c(k)$ as in Eq. \ref{eq:mod:V} and with the parameter estimates from Table \ref{tab:gg:S:d3}, the model $f_{S(C)}^{(16)}$  has the smallest AIC value and is accepted with a $p$-value of 0.87. 

\begin{table*}
\centering
\begin{tabular}{r|rr|rr}
$n$&$\hat{p}(n)$&s.\,e.&$\hat{p}_o(n)$&s.\,e.\\
\hline
4&0.188048&0.000055&0.004822&0.000022\\
5&0.258278&0.000062&0.036497&0.000059\\
6&0.225420&0.000059&0.102940&0.000096\\
7&0.156321&0.000051&0.171711&0.000119\\
8&0.091504&0.000041&0.201840&0.000127\\
9&0.046446&0.000030&0.182932&0.000122\\
10&0.020907&0.000020&0.134722&0.000108\\
11&0.008437&0.000013&0.083771&0.000087\\
12&0.003121&0.000008&0.045101&0.000066\\
13&0.001057&0.000005&0.021305&0.000046\\
14&0.000331&0.000003&0.009076&0.000030\\
15&0.000096&0.000001&0.003484&0.000019\\
16&0.000025&0.000001&0.001240&0.000011
\end{tabular}
\caption{Monte Carlo estimates of the probabilities ${p(n)=P(N(C)=n)}$ and ${p_o(n)=P(N(C_o)=n)}$ as well as corresponding standard errors (s.\,e.) in dimension ${d=3}$. Estimated probabilities for the higher side numbers not listed in the table: ${\hat{p}(17)=6.54\text{e-}6}$, ${\hat{p}(18)=1.42\text{e-}6}$, ${\hat{p}(19)=3.6\text{e-}7}$, ${\hat{p}(20)=6.0\text{e-}8}$, ${\hat{p}(21)=2.0\text{e-}8}$, ${\hat{p}_o(17)=3.928\text{e-}4}$, ${\hat{p}_o(18)=1.218\text{e-}4}$, ${\hat{p}_o(19)=3.54\text{e-}5}$, ${\hat{p}_o(20)=8.4\text{e-}6}$, ${\hat{p}_o(21)=2.0\text{e-}6}$, ${\hat{p}_o(22)=8.0\text{e-}7}$, ${\hat{p}_o(23)=1.0\text{e-}7}$; higher side numbers did not occur.}\label{tab:prob:d3}
\end{table*}

\begin{table*}
\centering
\begin{tabular}{r|r|rr|rr|rr|r}
$n$&$m(n)$&$\hat{\beta}_V(n)$&s.\,e.&$\hat{\nu}_V(n)$&s.\,e.&$\hat{\kappa}_V(n)$&s.\,e.&KS test\\
\hline
4&9402398&0.0337&0.0003&0.3356&0.0003&0.3001&0.0003&0.3216\\
5&12913883&0.0411&0.0004&0.6596&0.0006&0.2982&0.0003&0.1447\\
6&11271023&0.0530&0.0008&0.9927&0.0012&0.3015&0.0004&0.9080\\
7&7816070&0.0615&0.0014&1.3345&0.0023&0.3031&0.0006&0.6391\\
8&4575198&0.0712&0.0024&1.6726&0.0043&0.3052&0.0008&0.9973\\
9&2322319&0.0916&0.0048&1.9916&0.0078&0.3102&0.0013&0.7889\\
10&1045337&0.0877&0.0078&2.3529&0.0150&0.3082&0.0021&0.9695\\
11&421858&0.1052&0.0158&2.6727&0.0285&0.3116&0.0035&0.9769\\
12&156065&0.1267&0.0332&2.9893&0.0553&0.3154&0.0061&0.9488\\
13&52841&0.1006&0.0511&3.3898&0.1163&0.3094&0.0111&0.9728\\
14&16527&0.0300&0.0353&4.0218&0.2812&0.2856&0.0207&0.9031\\
15&4815&1.1957&1.4652&3.4529&0.3578&0.3753&0.0410&0.9801\\
16&1246&0.1815&0.6602&4.1400&1.0119&0.3181&0.0808&0.7087
\end{tabular}
\caption{ML estimates of the parameters $\beta$, $\nu$ and $\kappa$ of the generalized gamma distribution for ${V(C)\,\vert\,N(C)=n}$, the corresponding standard errors (s.\,e.) as well as the $p$-values of each a KS test.}\label{tab:gg:V:d3}
\end{table*}

\begin{table*}
\centering
\begin{tabular}{r|r|rr|rr|rr|r}
$n$&$m(n)$&$\hat{\beta}_S(n)$&s.\,e.&$\hat{\nu}_S(n)$&s.\,e.&$\hat{\kappa}_S(n)$&s.\,e.&KS test\\
\hline
4&9402398&1.1752&0.0056&0.5025&0.0004&0.4734&0.0005&0.4517\\ 
5&12913883&1.3287&0.0086&1.0014&0.0009&0.4727&0.0005&0.6279\\
6&11271023&1.5463&0.0140&1.5027&0.0018&0.4789&0.0007&0.9625\\
7&7816070&1.6524&0.0221&2.0159&0.0034&0.4807&0.0009&0.9908\\
8&4575198&1.7636&0.0361&2.5223&0.0063&0.4832&0.0013&0.9880\\
9&2322319&1.9747&0.0632&3.0043&0.0114&0.4888&0.0020&0.9730\\
10&1045337&1.8542&0.1015&3.5466&0.0221&0.4839&0.0032&0.9824\\
11&421858&2.0880&0.1924&4.0115&0.0417&0.4901&0.0054&0.9102\\
12&156065&2.2244&0.3621&4.5010&0.0816&0.4932&0.0094&0.9818\\
13&52841&1.8358&0.5859&5.1177&0.1727&0.4819&0.0170&0.9745\\
14&16527&0.8939&0.6502&6.0093&0.4102&0.4469&0.0316&0.7900\\
15&4815&8.7678&6.7297&5.1765&0.5261&0.5835&0.0627&0.9706\\
16&1246&1.5080&3.9464&6.5753&1.6742&0.4671&0.1232&0.8648
\end{tabular}
\caption{ML estimates of the parameters $\beta$, $\nu$ and $\kappa$ of the generalized gamma distribution for ${S(C)\,\vert\,N(C)=n}$, the corresponding standard errors (s.\,e.) as well as the $p$-values of each a KS test.}\label{tab:gg:S:d3}
\end{table*}

\subsubsection{The Crofton polyhedron}
Fitting only a single generalized gamma distribution to model the distribution of $V(C_o)$ (ML estimates ${\hat{\beta}_{V,o}=3.9313}$, ${\hat{\nu}_{V,o}=1.3834}$, ${\hat{\kappa}_{V,o}=0.3781}$) is again rejected by the KS test ($p$-value of $3.6\text{e-}6$). The same applies to $S(C_o)$ (ML estimates ${\hat{\beta}_{S,o}=20.4584}$, ${\hat{\nu}_{S,o}=2.2028}$, ${\hat{\kappa}_{S,o}=0.6028}$, $p$-value of $6.4\text{e-}18$). Therefore, again a more complex distribution model is required in each case.

Tables \ref{tab:gg:V0:d3} and \ref{tab:gg:S0:d3} provide the ML estimates ${(\hat{\beta}_{V,o}(n),\hat{\nu}_{V,o}(n),\hat{\kappa}_{V,o}(n))}$ and\linebreak ${(\hat{\beta}_{S,o}(n),\hat{\nu}_{S,o}(n),\hat{\kappa}_{S,o}(n))}$ of the parameters ${(\beta,\nu,\kappa)}$ of the generalized gamma distribution for ${V(C_o)\,\vert\,N(C_o)=n}$ and ${S(C_o)\,\vert\,N(C_o)=n}$, respectively, for $n=4,\ldots,19$, together with the related standard errors. As the $p$-values of the corresponding KS tests therein show, the conditional distributions ${V(C_o)\,\vert\,N(C_o)=n}$ or ${S(C_o)\,\vert\,N(C_o)=n}$ can be modeled sufficiently well by the corresponding generalized gamma distributions.

\medskip

Distribution of $V(C_o)$: Among the mixed models
\begin{equation}\label{eq:mod:V0}
\begin{aligned}
f_{V(C_o)}^{(k)}(x)=\sum_{n=4}^{k}\frac{\hat{p}_o(n)}{c_o(k)}\,f\bigg(\!x\,\bigg|\,\frac{\hat{\beta}_{V,o}(n)}{\lambda^3},\hat{\nu}_{V,o}(n),\hat{\kappa}_{V,o}(n)\!\!\bigg),&\\
\quad x>0,&
\end{aligned}
\end{equation}
with the parameter estimates from Tables \ref{tab:prob:d3} and \ref{tab:gg:V0:d3} and the sum $c_o(k)$ of all $\hat{p}_o(n)$ for $4\leq n\leq k$, the model $f_{V(C_o)}^{(19)}$ ($p$-value 0.53) is an appropriate approximation of the distribution of $V(C_o)$.

\medskip

Distribution of $S(C_o)$: Among the mixed models
\begin{equation}\label{eq:mod:S0}
\begin{aligned}
f_{S(C_o)}^{(k)}(x)=\sum_{n=4}^{k}\frac{\hat{p}_o(n)}{c_o(k)}\,f\bigg(\!x\,\bigg|\,\frac{\hat{\beta}_{S,o}(n)}{\lambda^2},\hat{\nu}_{S,o}(n),\hat{\kappa}_{S,o}(n)\!\!\bigg),&\\
\quad x>0,&
\end{aligned}
\end{equation}
with $c_o(k)$ as in Eq. \ref{eq:mod:V0} and with the parameter estimates from Table \ref{tab:gg:S0:d3} the model $f_{S(C_o)}^{(19)}$ ($p$-value 0.85) is an appropriate approximation of the distribution of $S(C_o)$.

\begin{table*}
\centering
\begin{tabular}{r|r|rr|rr|rr|r}
$n$&$m_o(n)$&$\hat{\beta}_{V,o}(n)$&s.\,e.&$\hat{\nu}_{V,o}(n)$&s.\,e.&$\hat{\kappa}_{V,o}(n)$&s.\,e.&KS test\\
\hline
4&241090&0.0456&0.0055&1.3031&0.0126&0.3081&0.0033&0.9995\\
5&1824825&0.0520&0.0027&1.6363&0.0065&0.3045&0.0013&0.9999\\
6&5147010&0.0672&0.0024&1.9665&0.0051&0.3077&0.0009&0.7401\\
7&8585535&0.0740&0.0023&2.3114&0.0051&0.3078&0.0007&0.9756\\
8&10091975&0.0865&0.0027&2.6442&0.0058&0.3100&0.0007&0.8139\\
9&9146575&0.0965&0.0034&2.9801&0.0072&0.3113&0.0008&0.8362\\
10&6736105&0.0973&0.0043&3.3385&0.0100&0.3107&0.0010&0.9684\\
11&4188555&0.1139&0.0067&3.6569&0.0146&0.3135&0.0013&0.9700\\
12&2255025&0.1216&0.0103&3.9960&0.0227&0.3143&0.0019&0.9854\\
13&1065260&0.1335&0.0173&4.3298&0.0371&0.3160&0.0028&0.9993\\
14&453800&0.1284&0.0272&4.6949&0.0645&0.3146&0.0045&0.5370\\
15&174195&0.2411&0.0784&4.8341&0.1064&0.3278&0.0075&0.9080\\
16&61985&0.2904&0.1627&5.1458&0.1960&0.3322&0.0131&0.8131\\
17&20098&0.1457&0.1667&5.6894&0.4098&0.3161&0.0234&0.9882\\
18&5985&0.1400&0.3106&6.1958&0.8539&0.3165&0.0447&0.7739\\
19&1717&2.3734&6.1853&5.1337&1.0828&0.3848&0.0842&0.8076
\end{tabular}
\caption{ML estimates of the parameters $\beta$, $\nu$ and $\kappa$ of the generalized gamma distribution for ${V(C_o)\,\vert\,N(C_o)=n}$, the corresponding standard errors (s.\,e.) as well as the $p$-values of each a KS test.}\label{tab:gg:V0:d3}
\end{table*}

\begin{table*}
\centering
\begin{tabular}{r|r|rr|rr|rr|r}
$n$&$m_o(n)$&$\hat{\beta}_{S,o}(n)$&s.\,e.&$\hat{\nu}_{S,o}(n)$&s.\,e.&$\hat{\kappa}_{S,o}(n)$&s.\,e.&KS test\\
\hline
4&241090&1.4107&0.1051&1.9756&0.0187&0.4841&0.0052&0.9381\\
5&1824825&1.4716&0.0482&2.5025&0.0098&0.4788&0.0021&0.9734\\
6&5147010&1.7137&0.0373&2.9973&0.0077&0.4850&0.0013&0.8530\\
7&8585535&1.7527&0.0333&3.5196&0.0076&0.4841&0.0011&0.7951\\
8&10091975&1.8965&0.0361&4.0133&0.0086&0.4872&0.0011&0.9825\\
9&9146575&1.9708&0.0428&4.5164&0.0108&0.4882&0.0012&0.8171\\
10&6736105&1.9388&0.0534&5.0493&0.0149&0.4866&0.0015&0.9636\\
11&4188555&2.0277&0.0752&5.5404&0.0217&0.4882&0.0020&0.9442\\
12&2255025&2.1249&0.1133&6.0282&0.0336&0.4902&0.0028&0.9844\\
13&1065260&2.1773&0.1777&6.5370&0.0552&0.4912&0.0043&0.9999\\
14&453800&2.1298&0.2833&7.0658&0.0954&0.4895&0.0068&0.6144\\
15&174195&3.2278&0.6554&7.2531&0.1565&0.5111&0.0114&0.9084\\
16&61985&3.4930&1.2316&7.7313&0.2895&0.5158&0.0200&0.5641\\
17&20098&2.2560&1.6203&8.5385&0.6039&0.4913&0.0358&0.9989\\
18&5985&1.1884&1.8992&9.8129&1.4078&0.4637&0.0681&0.8030\\
19&1717&8.0017&15.2450&8.2238&1.8205&0.5619&0.1286&0.9332
\end{tabular}
\caption{ML estimates of the parameters $\beta$, $\nu$ and $\kappa$ of the generalized gamma distribution for ${S(C_o)\,\vert\,N(C_o)=n}$, the corresponding standard errors (s.\,e.) as well as the $p$-values of each a KS test.}\label{tab:gg:S0:d3}
\end{table*}

\section{Discussion}
Based on extensive simulations, this paper shows that the distributions of the considered convex body functionals of the Poisson and the Crofton cell can be well modeled by a mixture of generalized gamma distributions; however, a single generalized gamma distribution is not sufficient. The chosen approach is motivated by the well-known decomposition of the distribution of the mean width according to the number of sides of these cells. The mixing weights are the probabilities for the occurrence of the individual side numbers. They were estimated from the simulations with great accuracy and are therefore available as part of the model parameters (Tables \ref{tab:prob:d2} and \ref{tab:prob:d3}).

Specifically, after selection through the AIC value, the following models can be used to model the distribution of the following variables: 
\begin{itemize}
\item area $A(C)$ of the Poisson polygon: $f_{A(C)}^{(9)}$, see Eq. \ref{eq:mod:A}
\item area $A(C_o)$ of the Crofton polygon: $f_{A(C_o)}^{(10)}$, see Eq. \ref{eq:mod:A0}
\item volume $V(C)$ of the Poisson polyhedron: $f_{V(C)}^{(16)}$, see Eq. \ref{eq:mod:V}
\item surface area $S(C)$ of the Poisson polyhedron: $f_{S(C)}^{(16)}$, see Eq. \ref{eq:mod:S}
\item volume $V(C_o)$ of the Crofton polyhedron: $f_{V(C_o)}^{(19)}$, see Eq. \ref{eq:mod:V0}
\item surface area $S(C_o)$ of the Crofton polyhedron: $f_{S(C_o)}^{(19)}$, see Eq. \ref{eq:mod:S0}
\end{itemize}

In addition to the fact that the (unconditional) distributions of these variables can be modeled as mixtures of generalized gamma distributions, it is equally remarkable that all (conditional) distributions of the mixture components can already be excellently modeled by special generalized gamma distributions (see the fit tests in Tables \ref{tab:gg:A:d2}, \ref{tab:gg:A0:d2}, \ref{tab:gg:V:d3}, \ref{tab:gg:S:d3}, \ref{tab:gg:V0:d3}, \ref{tab:gg:S0:d3}).

Table \ref{tab:gg:A:d2} shows that the values of $\hat{\nu}(n)$ for ${A(C)\,\vert\,N(C)=n}$ are approximately of the form ${(n-2)/2}$, i.\,e. roughly (and for $n>3$ even statistically indistinguishable) correspond to the first shape parameter of the generalized gamma distribution of ${B(C)^2\,\vert\,N(C)=n}$. Modeling the conditional distribution of ${A(C)\,\vert\,N(C)=n}$ by a generalized gamma distribution with the fixed parameter ${\nu=(n-2)/2}$ is therefore also possible (estimated parameter values are not reported here) and leads to suitable model candidates for $A(C)$ in a composition analogous to the mixture in Eq. \ref{eq:mod:A}. Similar observations and resulting variants with fixed parameters $\nu$ of the fitted generalized gamma distribution apply to $A(C_o)$ (${\nu=n/2}$), $V(C)$ (${\nu=(n-3)/3}$), $S(C)$ (${\nu=(n-3)/2}$), $V(C_o)$ ($\nu=n/3$) as well as $S(C_o)$ ($\nu=n/2$).

The models with the best overall AIC values are implemented in a Julia package \citep{Ballani2024b}, which can be used to quickly calculate the respective distribution functions and probability density functions, as well as to generate random numbers.

This is in a similar spirit to the work of \citet{Tanemura2003} and \citet{Shepilov2026}, in which comparable and similarly extensive simulation-based studies were made for the typical Poisson-Voronoi cell, another important model for random polyhedra \citep{SchneiderWeil2008,CSKM2013}.

\bibliographystyle{apalike-ejor}
\bibliography{poissoncell.bib}

@article{Ballani2024a,
  author = {Ballani, Felix},
  title = {{Data publication: Simulation results for standard statistics of Poisson and Crofton cells}},
  year  = {2024},
  journal = {Rodare},
  doi  = {10.14278/rodare.3290}
}

@article{Ballani2024b,
  author = {Ballani, Felix},
  title = {{DistributionModelsPHT: Julia package for distribution models for statistics of random cells of Poisson hyperplane tessellations (Version 0.4.0)}},
  year  = {2024},
  journal = {Rodare},
  doi  = {10.14278/rodare.3300}
}

@article{Ballani2024c,
  author = {Ballani, Felix},
  title = {{RandomCells: Julia package for the generation of random convex polytopes (Version 0.6.3)}},
  year  = {2024},
  journal = {Rodare},
  doi  = {10.14278/rodare.3231}
}

@article{BallaniBoogaart2014,
  author = {Ballani, F. and van den Boogaart, K.~G.},
  title = {{Weighted Poisson cells as models for random convex polytopes}},
  journal = {Methodology and Computing in Applied Probability},
  year = {2014},
  volume = {16},
  number = {2},
  pages = {369--384},
  doi = {10.1007/s11009-013-9342-y}
}

@article{BallaniStoyan2025,
  author = {Ballani, F. and Stoyan, D.},
  title = {Statistics for simulated assemblies of particles from mathematical models},
  journal = {Granular Matter},
  year = {2025},
  volume = {27},
  number = {46},
  doi = {10.1007/s10035-025-01519-6}
}

@article{BonnetEtAl2018,
  author = {Bonnet, G. and Calka, P. and Reitzner, M.},
  title = {{Cells with many facets in a Poisson hyperplane tessellation}},
  journal = {Advances in Mathematics},
  year = {2018},
  volume = {324},
  pages = {203--240},
  doi = {10.1016/j.aim.2017.11.016}
}

@article{Calka2001,
    author = {Calka, P.},
    title = {Mosa{\"i}ques poissoniennes de l'espace euclidien. {U}ne extension
	d'un r{\'e}sultat de {R}.~{E}.~{M}iles.},
    journal = {Comptes Rendus de l'Academie des Sciences, S{\'e}ries I},
    year = {2001},
    volume = {332},
    pages = {557--562},
    doi = {10.1016/S0764-4442(01)01885-7}
}

@article{Calka2003b,
    author = {Calka, P.},
    title = {Precise formulae for the distributions of the principal geometric characteristics of the typical cells of a two-dimensional {Poisson-Voronoi} tessellation and a {Poisson} line process},
    journal = {Advances in Applied Probability},
    year = {2003},
    volume = {35},
    pages = {551--562},
    doi = {10.1239/aap/1059486817}
}

@incollection{Calka2010,
    author = {Calka, P.},
    title = {Tessellations},
    booktitle = {New Perspectives in Stochastic Geometry},
    publisher = {Oxford University Press},
    address = {Oxford},
    year = {2010},
    editor = {Kendall, W.~S. and Molchanov, I.},
    pages = {145--169},
    doi = {10.1093/acprof:oso/9780199232574.003.0005}
}

@book{CSKM2013,
    title = {Stochastic {G}eometry and its {A}pplications},
    author = {Chiu, S.~N. and Stoyan, D. and Kendall, W.~S. and Mecke, J.},
    publisher = {J.~Wiley \& Sons},
    address = {Chichester},
    year = {2013},
    edition = {3rd},
    doi = {10.1002/9781118658222}
}

@article{CosterChermant2002,
    author = {Coster, M. and Chermant, J.-L.},
    title = {On a way to material models for ceramics},
    journal = {Journal of the European Ceramic Society},
    year = {2002},
    volume = {35},
    pages = {1191--1203},
    doi = {10.1016/S0955-2219(01)00455-1}
}

@article{CrainMiles1976,
    author = {Crain, I. and Miles, R.~E.}, 
    title = {Monte {C}arlo estimates of the distributions of the random polygons determined by random lines in the plane},
    journal = {Journal of Statistical Computation and Simulation},
    year = {1976},
    volume = {4},
    pages = {293--325},
    doi = {10.1080/00949657608810132}
}

@article{EscodaEtAl2015,
    author = {Escoda, J. and Jeulin, D. and Willot, F. and Toulemonde, C.},
    title = {Three-dimensional morphological modelling of concrete using multiscale {Poisson} polyhedra},
    journal = {Journal of Microscopy},
    year = {2015},
    volume = {258},
    pages = {31--48},
    number = {1},
    doi = {10.1111/jmi.12213}
}

@article{EscodaEtAl2016,
    author = {Escoda, J. and Willot, F. and Jeulin, D. and Sanahuja, J. and Toulemonde,
	C.},
    title = {Influence of the multiscale distribution of particles on elastic
	properties of concrete},
    journal = {International Journal of Engineering Science},
    year = {2016},
    volume = {98},
    pages = {60--71},
    doi = {10.1016/j.ijengsci.2015.07.010}
}

@article{George1987,
    author = {George, E.~I.},
    title = {Sampling random polygons},
    journal = {J. Appl. Prob.},
    year = {1987},
    volume = {24},
    pages = {557--573},
    doi = {10.2307/3214089}
}

@incollection{Grady2009,
    author = {Grady, D.},
    title = {{Dynamic Fragmentation of Solids}},
    booktitle = {{Shock Wave Science and Technology Reference Library: Volume 3, Solids II}},
    publisher = {Springer},
    address = {Berlin, Heidelberg},
    year = {2009},
    pages = {169--276},
    doi = {10.1007/978-3-540-77080-0_4}
}

@book{Grady2017,
    author = {Grady, D.},
    title = {{Physics of Shock and Impact. Volume 1: Fundamentals and dynamic failure}},
    publisher = {{IOP Publishing}},
    address = {{Bristol, UK}},
    year = {2017},
    doi = {10.1088/978-0-7503-1254-7}
}

@book{HugSchneider2024,
    author = {Hug, D. and Schneider, R.},
    title = {{Poisson Hyperplane Tessellations}},
    publisher = {Springer Monographs in Mathematics},
    address = {Berlin},
    year = {2024},
    doi = {10.1007/978-3-031-54109-9}
}

@book{Jeulin2021,
  author = {Jeulin, D.},
  title = {Morphological Models of Random Structures},
  publisher = {Springer},
  address = {Cham},
  year = {2021},
  doi = {10.1007/978-3-030-75452-5}
}

@incollection{KlattEtAl2017,
    author = {Klatt, M.~A. and Last, G. and Mecke, K. and Redenbach, C. and Schaller,
	F.~M. and Schr\"oder-Turk, G.~E.},
    title = {Cell shape analysis of random tessellations based on {M}inkowski
	tensors},
    booktitle = {Tensor Valuations and Their Applications in Stochastic Geometry and
	Imaging},
    publisher = {Springer International Publishing},
    address = {Cham},
    year = {2017},
    editor = {Jensen, E.~B.~V. and Kiderlen, M.},
    pages = {385--421},
    doi = {10.1007/978-3-319-51951-7_13}
}

@book{Lantuejoul2002,
    title = {Geostatistical {S}imulation},
    author = {Lantu{\'e}joul, C.},  
    publisher = {Springer, Berlin},
    address = {Berlin},
    year = {2002},
    doi = {10.1007/978-3-662-04808-5}
}

@article{LantuejoulEtAl2011,
    author = {Lantu{\'e}joul, C. and Bacro, J.-N. and Bel, L.},
    title = {Storm processes and stochastic geometry},
    journal = {Extremes},
    year = {2011},
    volume = {14},
    pages = {413--428},
    doi = {10.1007/s10687-010-0121-7}
}

@article{MichelParoux2007,
    author = {Michel, J. and Paroux, K.},
    title = {Empirical polygon simulation and central limit theorems for the homogeneous {Poisson} line process},
    journal = {Methodol. Comput. Appl. Probab.},
    year = {2007},
    volume = {9},
    pages = {541--556},
    doi = {10.1007/s11009-006-9009-z}
}

@article{Miles1964a,
    author = {Miles, R.~E.},
    title = {Random polygons determined by random lines in a plane},
    journal = {Proc. Nat. Acad. Sci. (USA)},
    year = {1964},
    volume = {52},
    pages = {901--907},
    doi = {10.1073/pnas.52.4.901}
}

@article{Miles1971,
    author = {Miles, R.~E.},
    title = {Poisson flats in {E}uclidean spaces. {P}art {II}: {H}omogeneous {P}oisson
	flats and the complementary theorem.},
    journal = {Advances in Applied Probability},
    year = {1971},
    volume = {3},
    pages = {1--43},
    doi = {10.2307/1426328}
}

@incollection{Miles1974,
    author = {Miles, R.~E.},
    title = {A synopsis of `{P}oisson flats in {E}uclidean spaces'},
    booktitle = {Stochastic Geometry},
    publisher = {John Wiley \& Sons},
    address = {London},
    year = {1974},
    editor = {Harding, E.~F. and Kendall, D.~G.},
    pages = {202--227}
}

@article{MollerZuyev1996,
    author = {M{\o}ller, J. and Zuyev, S.},
    title = {Gamma-type results and other related properties of {P}oisson processes},
    journal = {Advances in Applied Probability},
    year = {1996},
    volume = {28},
    pages = {662--673},
    doi = {10.2307/1428175}
}

@techreport{MottLinfoot1943,
    author = {Mott, N.~F. and Linfoot, E.~H.},
    title = {A {T}heory of {F}ragmentation},
    institution = {United Kingdom Ministry of Supply},
    year = {1943},
    number = {AC3348},
    month = {January},
    doi = {10.1007/978-3-540-27145-1_9}
}

@article{NagelWeiss2003,
    author = {Nagel, W. and Weiss, V.},
    title = {Limits of sequences of stationary planar tessellations},
    journal = {Advances in Applied Probability (SGSA)},
    year = {2003},
    volume = {35},
    pages = {123--138},
    doi = {10.1239/aap/1046366102}
}

@article{NagelWeiss2005,
    author = {Nagel, W. and Weiss, V.},
    title = {Crack {STIT} tessellations: characterization of stationary random
	tessellations stable with respect to iteration},
    journal = {Advances in Applied Probability},
    year = {2005},
    volume = {37},
    pages = {859--883},
    doi = {10.1239/aap/1134587744}
}

@incollection{QuenecEtAl1994,
    author = {Quenec'h, J.~L. and Chermant, J.~L. and Coster, M. and Jeulin, D.},
    title = {Liquid phase sintered materials modelling by random closed sets},
    booktitle = {Mathematical Morphology and its Applications to Image Processing},
    publisher = {Kluwer Academic Pub.},
    address = {Dordrecht},
    year = {1994},
    editor = {J. Serra and P. Soille},
    pages = {225--232},
    doi = {10.1007/978-94-011-1040-2_29}
}

@article{RebbahEtAl2019,
    author = {Rebbah, S. and Nicol, F. and Puechmorel, S.},
    title = {The geometry of the generalized gamma manifold and an application to medical imaging},
    journal = {Mathematics},
    year = {2019},
    volume = {7},
    pages = {674},
    doi = {10.3390/math7080674}
}

@book{SchneiderWeil2008,
    title = {{Stochastic and Integral Geometry}},
    author = {Schneider, R. and Weil, W.},
    publisher = {Springer},
    address = {Berlin},
    year = {2008},
    doi = {10.1007/978-3-540-78859-1}
}

@book{Serra1982,
    title = {{Image Analysis and Mathematical Morphology}},
    author = {Serra, J.},
    publisher = {{Academic Press}},
    address = {London},
    year = {1982}
}

@article{Shepilov2026,
  author = {Shepilov, M.},
  title = {On analytical approximation of cell volume distribution in three-dimensional {Poisson-Voronoi} tessellation},
  journal = {Annals of Mathematics and Physics},
  year = {2026},
  volume = {9},
  number = {4},
  pages = {231--239},
  doi = {10.17352/amp.000200}
}

@article{Stacy1962,
    author = {Stacy, E.~W.},
    title = {A generalization of the gamma distribution},
    journal = {Annals of Mathematical Statistics},
    year = {1962},
    volume = {33},
    number = {3},
    pages = {1187--1192},
    doi = {10.1214/aoms/1177704481}
}

@article{Tanemura2003,
  author = {Tanemura, M.},
  title = {Statistical distributions of {P}oisson {V}oronoi cells in two and three dimensions},
  journal = {Forma},
  year = {2003},
  volume = {18},
  pages = {221--247}
}

@article{Tanner1983a,
    author = {Tanner, J.~C.},
    title = {The proportion of quadrilaterals formed by random lines in a plane},
    journal = {Journal of Applied Probability},
    year = {1983},
    volume = {20},
    pages = {400--404},
    doi = {10.2307/3213813}
}

\end{document}